# ANALYSIS ON THE PROJECTIVE OCTAGASKET

YIRAN MAO, ROBERT S. STRICHARTZ, LEVENTE SZABO, AND WING HONG WONG

ABSTRACT. The existence of a self similar Laplacian on the Projective Octagasket, a non-finitely ramified fractal is only conjectured. We present experimental results using a cell approximation technique originally given by Kusuoka and Zhou. A rigorous recursive algorithm for the discrete Laplacian is given. Further, the spectrum and eigenfunctions of the Laplacian together with its symmetries are categorized and utilized in the construction of solutions to the heat equation.

## Contents

## 1. INTRODUCTION

A general method due to Kigami[Ki], for constructing Laplacians on self-similar fractals involves realizing the fractal as a limit of finite graphs, and attempt to define the fractal Laplacian as a limit of graph Laplacians. This is usually successful in the case of *postcritically finite* (p.c.f) fractals which is a loose notion of *finitely ramified* fractals. Examples include the Sierpinski gasket, and various other gasket fractals based on regular polygons. These fractals are connected but just barely

| Yiran Mao | Robert S. Strichartz | Levente Szabo | Wing Hong Wong |
|---|---|---|---|
| Math Dept | Math Dept | Math Dept | Math Dept |
| Gettysburg College | Cornell University | UNC Asheville | The Chinese University of Hong Kong |
| Gettysburg PA 17325 | Ithaca NY 14853 | Asheville NC 28804 | Hong Kong |
| maoyi01@gettysburg.edu | str@math.cornell.edu | lszabo@unca.edu | 1155064377@link.cuhk.edu.hk |

connected, and the removal of a finite number of points disconnects them. In contrast, the octagasket (Figure 1) built from a regular octagon, is clearly not p.c.f., and only recently [A] has a rigorous construction of a Laplacian been attempted using probabilistic methods, although several earlier experimental works have studied spectral properties of the octagasket Laplacian.

One source of ambiguity in this area involves how to define and deal with a notion of *boundary*. One way to attempt to deal with this problem is to replace the original fractal with a fractal without boundary by identifying points. In a classical setting one would replace a square by a torus ( or projective space or Klein bottle) by identifying points on the boundary of the square, leading to a manifold without boundary. A number of works[BKS,MOS] have considered such modifications of classical fractals. In particular, many classical fractals appear to have a countable set of "boundary" pieces, and the *projective* version consists of identifying antipodal points on each of the pieces. In this paper we study the *projective octagasket* (POG) obtained by applying this procedure to the classical octagasket.

In order to study a possible Laplacian on POG we use an idea due to Kusuoka and Zhou [KZ] to work with *cell graph* approximations. At each approximation level $m$ we subdivide the POG into $8^m$ *m-cells* in the obvious way. The cell graph at level $m$ is made up of vertices equal to the set of $m$-cells with an edge connecting two cells that have an edge in common; in particular there may be more than one edge connecting a pair of vertices, and each vertex has exactly 8 edges connected to it. given a continuous function on the POG, we obtain a function on the cell graph by assigning to each $m$-cell the average value on the 8 vertices of the associated octagon. Let $\triangle_m f$ be the graph Laplacian on the m-cell graph. It is reasonable to expect that there exist renormalization constants $r_m$ (perhaps $r_m = r^m$ for some $r$) such that $r_m \triangle_m f$ converges to an operator $\triangle_m f$ on POG for reasonably well behaved functions $f$.

The main purpose of this paper is to present experimental evidence for the validity of these expectations, to estimate the renormalization constants $r_m$, to compute approximately the spectra of the graph Laplacians and to deduce from this data properties of the spectrum of the limiting POG Laplacian. To this end we give the construction of the POG in section 2 and the definition of the Laplacian $\triangle_m$ on the $m$-cell graph in section 3. The heart of the matter is the identification algorithm described in section 4 that gives a precise recursive labelling scheme for the $m$-cell graph. With this in hand it is hand it is fairly straightforward to compute the eigenvalues and eigenfunctions in section 5. We report some this data here, but refer the reader to website [W] for the programs used in the computations and the complete data sets. In section 6 we show graphs of some of the eigenfunctions of multiplicity one, as evidence of the convergence of the eigenfunctions of level m, as $m \to \infty$. We also discuss the behavior of the graphs of the eigenvalue counting function and Weyl ratio as $m$ varies. In sections 7 and 8 we discuss applications to solutions of the heat equation and wave equation on POG. Finally, section 9 is devoted to some basic questions on the geometry of the $m$-cell graphs. There is a natural metric in terms of shortest paths between $m$-cells. Specifically, what is the diameter of the graph? How many $m$-cells are contained in the ball of given radius about an $m$-cell? We formulate some conjectured answers and give proofs of some weaker statements to the effect that the number of cells is on the order of $n^3$ where $n$ is the radius.

We conclude the paper in section 10 with a discussion of the extend to which our experimental evidence confirm our initial expectations.

## 2. Construction of Projective Octagasket

To define the projective octagasket, we first need to describe the octagasket. The octagasket is defined to be the unique nonempty compact set $OG$ satisfying the following identity with contraction mappings $\{F_n\}$:

$$OG = \bigcup_{n=0}^{7} F_n(OG)$$

where

$$F_n(x,y) = \frac{1}{2+\sqrt{2}}(x + (1+\sqrt{2})\cos\frac{5-2n\pi}{8}, y + (1+\sqrt{2})\sin\frac{5-2n\pi}{8})$$

for all $x, y \in \mathbb{R}$ and all $n = 0, 1, ..., 7$.

Of course, the octagasket can be obtained through finite graph approximations. Let $\Gamma_0$ be a regular octagon centered at the origin with radius 1, $V_0$ be its vertices and $R_0$ be the origin. Then, we recursively let

$$\Gamma_m = \bigcup_{n=0}^{7} F_n(\Gamma_{m-1}) \quad \forall m \in \mathbb{N}$$

and

$$V_m = \bigcup_{n=0}^{7} F_n(V_{m-1}) \quad \forall m \in \mathbb{N}$$

One can easily see that $V_* := \bigcup_{n=0}^{\infty} V_n$ is dense in $OG$ and $\bigcap_{n=0}^{\infty} \Gamma_n = \overline{V_*} = OG$. The following shows figures 2.0 - 2.5, namely $\Gamma_0, \Gamma_1, \Gamma_2, \Gamma_3, \Gamma_4, \Gamma_5$ respectively.

**Figure 2. Levels 0-5**

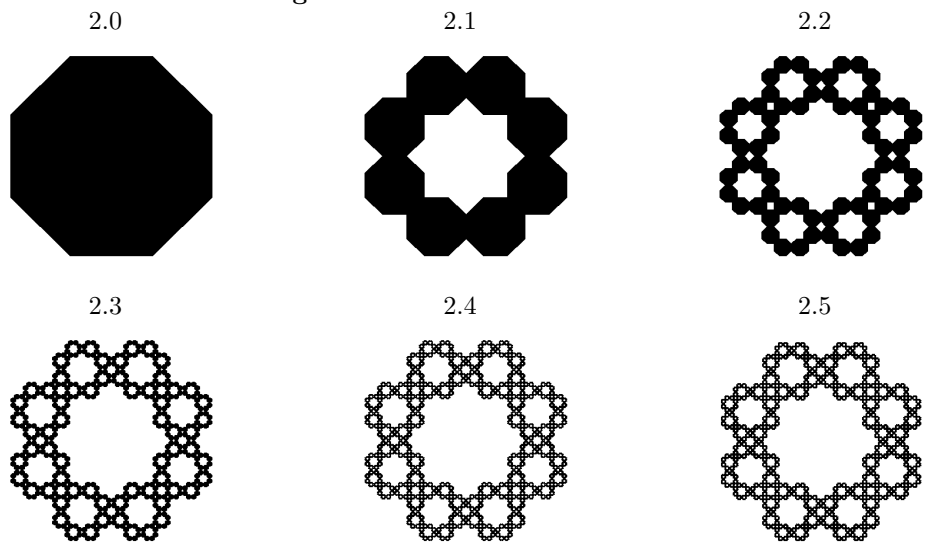

2.0  2.1  2.2

2.3  2.4  2.5

$K$ is called an $m$-cell of $OG$ if $K = F_{j_1}Fj_2...Fj_m(OG)$ for some $j_1, j_2, ..., j_m = 0, 1, ..., 7$. $L$ is called a pseudo-$m$-cell if $L$ is not an $m$-cell but there exists a rigid motion $\phi$ such that $\phi(L)$ is an $m$-cell.

Intuitively, the projective octagasket is given by identifying the boundaries of the octagasket antipodally. This process is depicted in Figure 2.7. Two points are identified if they belong to the inner boundaries of the same $m$-cell and they are antipodal with respect to the centre of that $m$-cell; or they belong to the outer boundaries of the 0-cell and they are antipodal with respect to the origin; or they both belong to the inner boundary of the same pseudo-$m$-cell and they are antipodal with respect to the centre of that pseudo-$m$-cell. In the following graph, the red, blue and green lines are corresponding to some gluing of inner boundaries of an $m$-cell, outer boundaries of the 0-cell and inner boundaries of pseudo-$m$-cells on $\Gamma_3$, note the tiny squares are glued since they will become a pseudo-3-cell eventually.

**Figure 2.7 Antipodal Identifications**

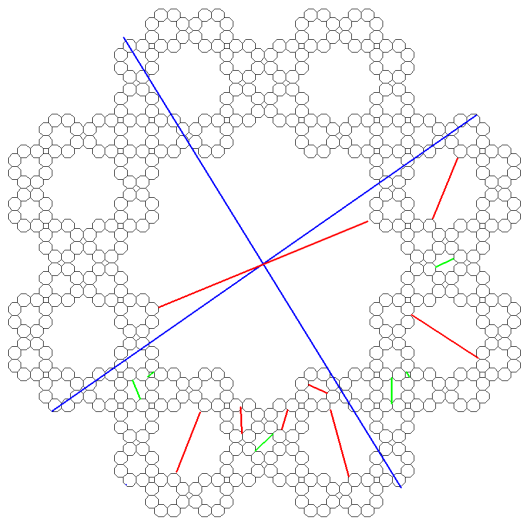

To define the projective octagasket, we first let $E_*$ denote the collection of all centres of all the pseudo-cells and $R_*$ is the collection of all the centres of all the cells.

To be precise, we let $\rho = \frac{1}{2+\sqrt{2}}$, $q_n$ be rotation by $\frac{n\pi}{4}$ anti-clock-wisely centered at the origin, $R_0$ be the origin,

$$P = \{(0, (1-\rho)\cos\frac{\pi}{8} + \sum_{i=1}^{n} \epsilon_i \rho^i (1-\rho)\sin\frac{\pi}{8}) \in \mathbb{R}^2 | \epsilon_i = \pm 1, n = 0, 1, 2, ...\}$$

$$E_0 = \bigcup_{n=0}^{7} F_n(P)$$

$$R_m = \bigcup_{n=0}^{7} F_n(R_{m-1}) \quad \forall m \in \mathbb{N}$$

$$E_m = \bigcup_{n=0}^{7} F_n(E_{m-1}) \quad \forall m \in \mathbb{N}$$

where $P$ is the collection of all centres of all the pseudo-cells between the 1-cells 0 and 1, $E_0$ is the collection of all centres of all the pseudo-cells between the 1-cells. Then, $R_* = \bigcup_{n=0}^{\infty} R_n$ and $E_* = \bigcup_{n=0}^{\infty} E_n$. The red points in Figure 2.8 below denotes the location of the points in $P$:

**Figure 2.8 P: Pseudo-cell centers**

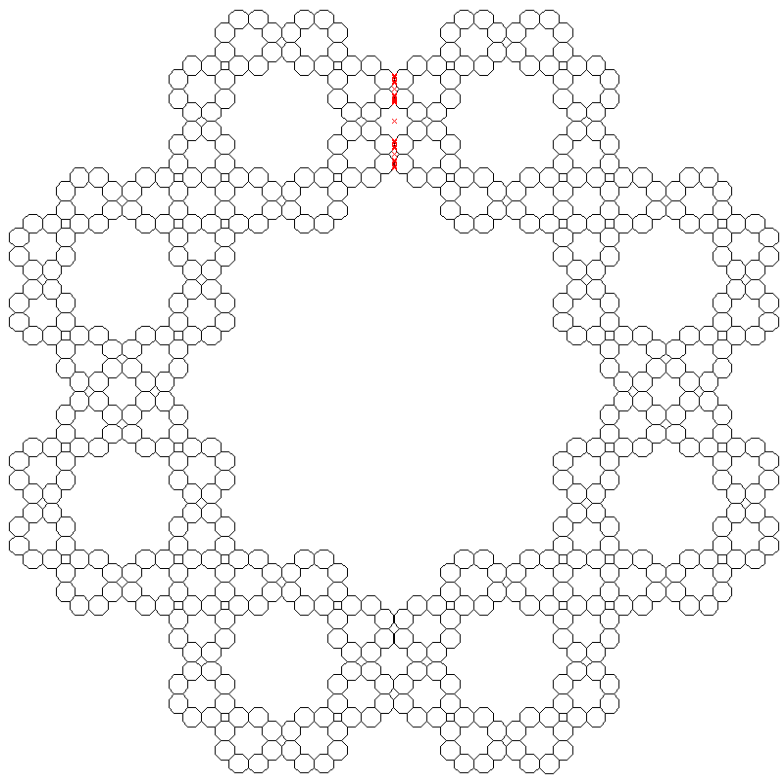

Now, define an equivalence relation on $OG$ by $x \sim y$ if $l_{xy} \cap (R_* \cup E_*) \neq \emptyset$ and $l_{xy} \cap OG = \{x, y\}$, where $l_{xy}$ is a straight line in the $\mathbb{R}^2$ with end points $x$ and $y$; or there exits $l_{xy}$ such that $l_{xy} \cap R_0 \neq \emptyset$ with any extension $\widetilde{l_{xy}}$ of $l_{xy}$ satisfying $\widetilde{l_{xy}} \cap OG = l_{xy} \cap OG$. Then, we define the projective octagasket $POG$ to be $OG/\sim$.

*Remark* 2.1. $l_{xy} \cap R_* \neq \emptyset, l_{xy} \cap E_* \neq \emptyset$ and $\widetilde{l_{xy}} \cap OG = l_{xy} \cap OG$ correspond to the identification of inner boundaries of $m$-cell, the identification of inner boundaries of pseudo-$m$-cell and the identification of outer boundary of 0-cell respectively.

## 3. Laplacian

The Laplacian plays a central role in constructing an analytic theory on a given fractal. In general a fractal is given as a limit of graph approximations. As such a discrete graph Laplacian is used to model the continuous Laplacian. When a renormalisation factor is found then the Laplacian at a point can be given as the limit of its discrete graph approximations. Generally the graph Laplacian is given

as
$$-\triangle_m u(x) = \sum(u(y) - u(x))$$
over all vertices $y$ that neighbor $x$ in the $m$ level graph approximation. We will use a form of cell graph approximation given by Kusuoka and Zhou in [2]. Here the function value on an $m$ cell $A_i^m$ is given as the average value on its vertices. Now given some $u : V \to \mathbb{R}$ we can evaluate functions on cells and therefore we can construct a pointwise Laplacian as
$$-\triangle_m u(A_x^m) = \sum(u(A_y^m) - u(A_x^m))$$
Where we sum over the 8 neighboring cells $A_j^m$. Every $m$ cell has 8 edges, each bordering another $m$ cell. Thus in total every cell will have 8 adjacent cells counting multiplicities. The Laplacian matrix is an $8^m \times 8^m$ matrix; for level 1 the matrix is given in Figure 3.0.

**Figure 3.0 Laplacian Matrix: Level 1**

$$\begin{bmatrix} 8 & -1 & 0 & 0 & -6 & 0 & 0 & -1 \\ -1 & 8 & -1 & 0 & 0 & -6 & 0 & 0 \\ 0 & -1 & 8 & -1 & 0 & 0 & -6 & 0 \\ 0 & 0 & -1 & 8 & -1 & 0 & 0 & -6 \\ -6 & 0 & 0 & -1 & 8 & -1 & 0 & 0 \\ 0 & -6 & 0 & 0 & -1 & 8 & -1 & 0 \\ 0 & 0 & -6 & 0 & 0 & -1 & 8 & -1 \\ -1 & 0 & 0 & -6 & 0 & 0 & -1 & 8 \end{bmatrix}$$

Similar to what has been done on the Sierpinski gasket and other fractals, after constructing the Laplacian $\triangle_m$ for $\Gamma_m$, it is natural to construct the Laplacian for the projective octagasket by letting
$$\triangle u(x) = \lim_{m \to \infty} r^{-m} \triangle (A_i^m)$$
for some renormalization factor $r$, where $x \in A_i^m$ for all m and where $r$ is determined later.

## 4. Identification Algorithm

We now turn to the labelling algorithm which lies at the core of constructing the graph approximation. We wish to construct the Laplacian for a given $m$ level cell graph approximation. To do this we need a labelling scheme and then a recursive algorithm which allows us to find the number of edge connections between neighboring cells.

Recall that $\Gamma_m$ denotes the $m$ level cell graph approximation of the POG. For $\Gamma_1$ we label the 1 cells clockwise $0-7$ starting at the top left. For $\Gamma_2$ we take the cell labelled 0 on the previous level and label its 2 cells $0-7$ as before. Now for the cell labelled 1 we rotate our labelling scheme clockwise by $\pi/4$, thus labelled $0-7$ clockwise starting on the top right. This process is repeated over $\Gamma_2$.

For $\Gamma_m$ simply take the labelling on $\Gamma_{m-1}$ and to label cell number $k$ simply rotate the $0-7$ clockwise labelling scheme by $k\pi/4$

To be rigorous, $w$ is an $m$-cell in $\Gamma_m$ if $w = F_{j_1}F_{j_2}...F_{j_m}(\Gamma_0)$ for some $j_1, j_2, ..., j_m$ elements of $0, 1, ..., 7$ and we write $w = j_1(j_2-j_1)...(j_m-j_{m-1})$ and $j_1(j_2-j_1)...(j_m-j_{m-1})$ is called the address of $w$.

In addition, $\nu$ is called a pseudo-$m$-cell in $\Gamma_{m+k}$ if there exists some rigid motion $\phi$ such that $\phi(\nu)$ is an $m$-cell but $\nu$ is not an $m$-cell for $k > 0$; or $\nu$ is a square for $k = 0$. The pseudo-cells on $\Gamma_m$ are essentially the boundary cells between the $m$ cells.

**Through out this paper, any arithmetic on the address will be conducted in $\mathbb{Z}_8$.**

Figures 4.0, 4.1, 4.2 show the addressing on $\Gamma_1, \Gamma_2$ and $\Gamma_3$ respectively:

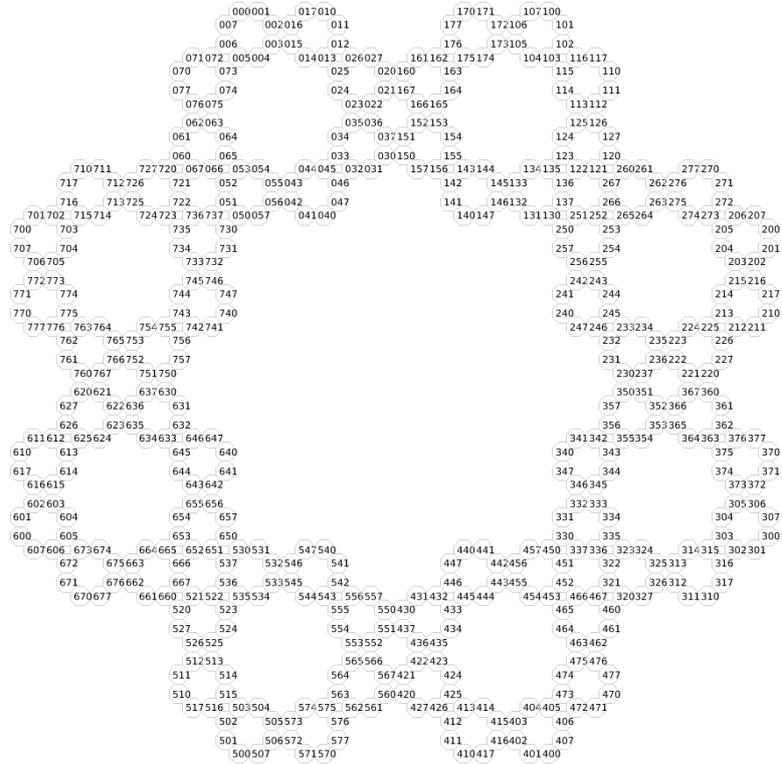

**Figure 4.2 Address Level 3**

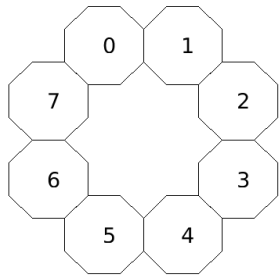

**Figure 4.0 Address Level 1**

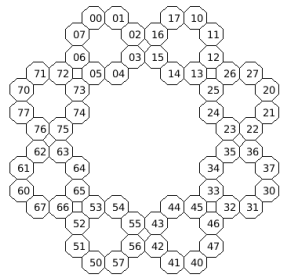

**Figure 4.1 Address Level 2**

The number of edges shared between cells is crucial to the Laplacian matrix of our approximation. Let $A_i^m$ be an m-cell. We initialize a $8^m \times 8^m$ matrix $L$ where row $i$ will store the value of $-\triangle A_i^m$ in terms of the other $m$ cells. Thus $L(i,i) = 8$ for all $i$ and $L(i,j)$ will be negative of the number of edges between the two cells $i$ and $j$. Recall that our labelling algorithm gives an $m$ cell as a $m$ tuple of the references. Thus each cell can be represented in base 8. Now in order to reference a cell in the Laplacian matrix we simply translate it to base 10. Thus some $(a,b,c)$ will correspond to row $64a + 8b + c$ in the matrix $L$.

The identification for $\Gamma_1$ is simple, for cell number $i$ there is an edge shared between both $(i-1)$ and $(i+1)$. There are 6 edges connecting it to $(i+4)$, 2 lie on the inner boundary and 4 lie on the outer boundary. This allows one to construct the Laplacian matrix easily.

To understand the structure of $\Gamma_{m+1}$ for $m > 0$, we will need to understand when two cells share an edge. We have the following algorithm to create $\Gamma_{m+1}$ from $\Gamma_m$:

(1) Remove the outer boundary identification of $\Gamma_m$.
(2) Make 8 copies of $\Gamma_m$.
(3) Identify the adjacent edges of the $\Gamma_m$ and the inner boundary of the pseudo-cell antipodally.
(4) Identify the inner boundary of the 0-cell antipodally
(5) Identify the outer boundary of the 0-cell antipodally.

(1) **Remove outer boundary identification**
    For $m > 1$, we simply omit step (5) for the previous construction. For $m = 1$, it is easy.
(2) **Create 8 copies of $\Gamma_m$**
    Simply take the disjoint union of 8 copies of $\Gamma_m$, label them according to the previous labelling scheme.
(3) **Adjacent Edges and Pseudo-cell Identifications**
    This is the most difficult part of the identification algorithm. We show the explicit construction for level 2,3,4 and a recursive generalization for higher levels. On $\Gamma_2$ one can see the small pseudo cells, they are the inner boundary of four adjacent cells. Notice that for any two adjacent 2 cells these pseudo-cell boundaries always lie between $X2, X3$ and $(X+1)5, (X+1)6$. We will start our explanation here. Notice that for each subsequent level the number of pseudo-cells doubles. Also notice that if any cell has $n$ edges along the boundary then on the next level those $n$ edges will be replaced by $n + 1$ cells.

**Level 2**

Initially we take the adjacency gluing of $\Gamma_2$, this is where the four cells are glued along a square antipodally. Given an arbitrary 1 cell denoted $X$ from 0-7 we construct the gluing below. Here every edge denotes 1 edge connection between cells. In Figure 4.3 we show the pseudo-cells from $\Gamma_2$. In Figure 4.4 we give the labelling mechanism across the pseudo-cell for any arbitrary adjacent 1 cells in $\Gamma_2$.

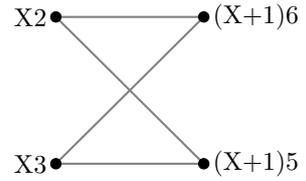

**Figure 4.4 Level 2 Identifications**

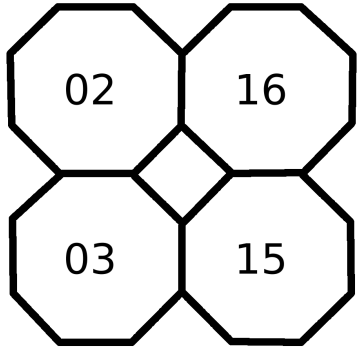

**Figure 4.3 Level 2 Pseudo-cell**

**Level 3**

Recall that for a given cell with $n$ edges along the boundary on the next level that cell will turn into $n+1$ finer cells. To construct the next level we observe how many edge connections each cell has. In $\Gamma_2$ each has 2, so we take 2+1 copies of each cell and append to it the repeating string $\overline{01267}$ where we start at 0 and continue with this pattern throughout the cell. When changing over to a next cell a 0 will follow a 0 and a 6 follows a 2. This process is for labelling the left side, the right side is identical except with the order reversed. The gluing process is copied twice from $\Gamma_2$ here we have $X20$ together with $X21$ paired with $X60$ and $X67$ in the cross formation as well as $X37$ and $X30$ paired with $X51$ and $X50$ as shown in black. The central 4 cells are connected antipodally, as shown in red. All of these central edge connections correspond to 2 edges between cells. In Figure 4.5 we show the pseudo-cells for $\Gamma_3$. In Figure 4.6 we give the labelling mechanism across the pseudo-cell for any arbitrary adjacent 1 cells in $\Gamma_3$.

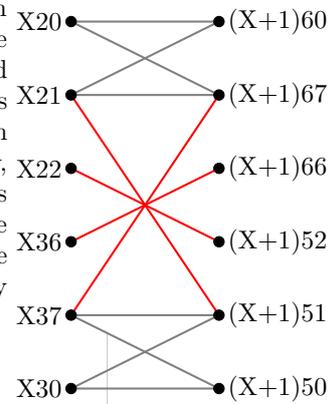

Figure 4.6 Level 3 Identifications

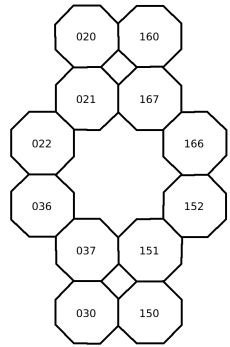

Figure 4.5 Level 3 Pseudo-cell

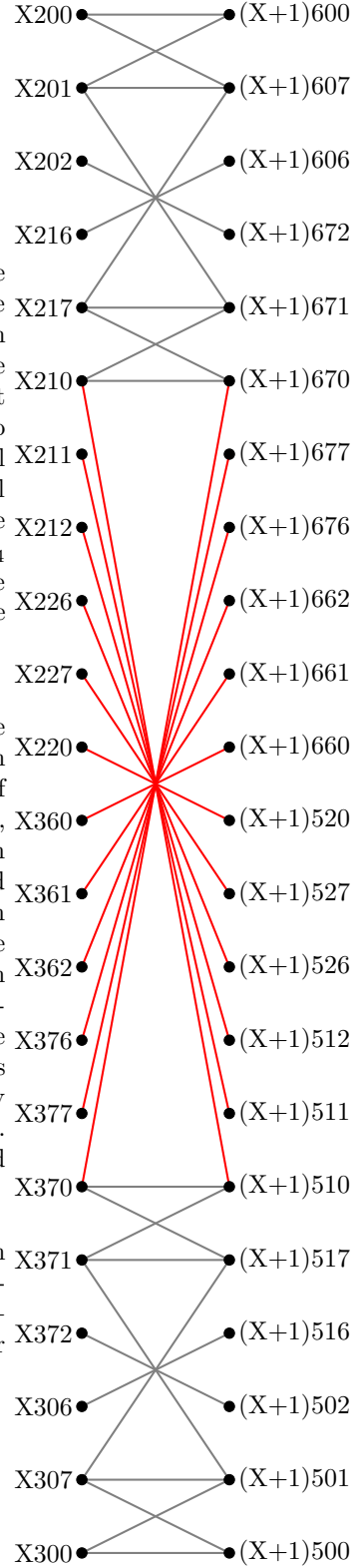

**Figure 4.7 Level 4 Identification**

**Level 4**
For $\Gamma_4$ we again take $n + 1$ copies of each cell where $n$ is the number of edge connections of the cell in the step (3) of previous level. So a cell like $X21$ from before has two edges each of size 1 and an edge of size 2, thus there will be 5 copies of $X21Y$ on the next level. The gluing structure of $\Gamma_3$ is copied over to the top and bottom of the $\Gamma_4$ boundary. The central 12 cells are glued antipodally. On the previous level we had all central edge connections equal to 2, we had a sequence $2, 2, 2, 2$ of edge connections. On $\Gamma_4$ we replace every 2 in this sequence with the sequence $242$, thus getting $2, 4, 2, 2, 4, 2, 2, 4, 2, 2, 4, 2$ as the edge connections across the center.

**Level m**
To give a recursive argument we assume we have the adjacency labelling along $\Gamma_{m-1}$ now each cell has an edge connection $n$ associated with it in the step (3) of previous level. For $\Gamma_m$ take $n + 1$ copies of each cell, list everything in the same order as before and begin labelling with $\overline{01267}$. Whenever a new cell is reached we restart the labelling at 0 if the previous one ends in 0 and a 6 if the previous one ends in 2. Now copy the exact gluing structure of $\Gamma_{m-1}$ on the top and bottom parts of $\Gamma_m$. The central portion is connected antipodally as before and has a edge connection sequence associated with it. This sequence from $\Gamma_{m-1}$ consists of 2's and 4's, whenever there is a 2 we replace by $2, 4, 2$ whenever there is a 4 we replace by $2, 4, 4, 4, 2$. This gives the new edge strength sequence for $\Gamma_m$ and we are finished.

(4) **Inner Boundary identifications**
If an $(m + 1)$-cell $w = X_1X_2X_3...X_{m+1} \in \Gamma_{m+1}$ with $X_2 = 3, 4$ or 5 shares $n$ edges with other cells before (4), then $w$ shares $(8 - n)$ edges with $(X_1 + 4)X_2X_3...X_{m+1}$ after (4). The connectivity of other cells remains unchanged.

## Figure 4.8 Level 3 Outer Boundary

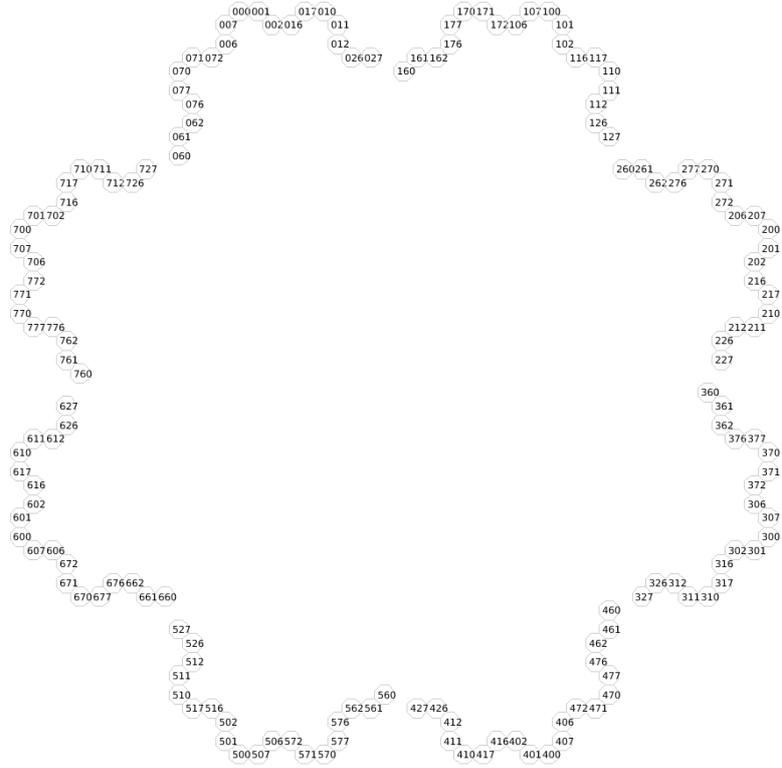

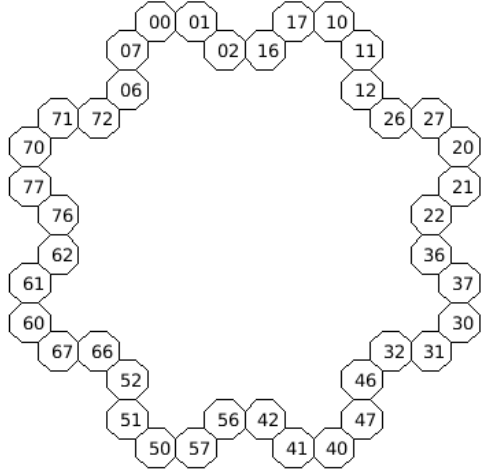

**Figure 4.9 Level 2 Outer Boundary**

(5) **Outer Boundary identifications**

If an $(m+1)$-cell $w = X_1 X_2 X_3 ... X_{m+1} \in \Gamma_{m+1}$ shares $n$ edges with other cells before (5), then $w$ shares $(8 - n)$ edges with $(X_1 + 4) X_2 X_3 ... X_{m+1}$ after (5).

We show only the outer boundary for level 2 and 3 in figures 4.9 and 4.8 respectively. If one attempts to connect these cells antipodally it is apparent that one is doing something equivalent to the algebraic identification listed above.

## 5. Eigenvalue and Eigenfunction

Recall that we will attempt to define the Laplacian $-\triangle u(x) = -\lim_{m\to\infty} r^{-m} \triangle_m u(A_i^m).$, where $x \in A_i^m$ for all $m$. In order to guess the renormalization factor $r$, we first compute the eigenvalue and eigenfunction of $-\triangle_m$.

In Table 5.0 a portion of the eigenvalues for levels 1, 2, 3, 4, 5 are listed. Observe that the multiplicity of every eigenvalue does not exceed 2.

**Table 5.0 Eigenvalues**

| Level1 Eigenvalue | | Level2 Eigenvalue | | Level3 Eigenvalue | | Level4 Eigenvalue | | Level5 Eigenvalue | |
|---|---|---|---|---|---|---|---|---|---|
| 0 | 1 | 0 | 1 | 0 | 1 | 0 | 1 | 0 | 1 |
| 2 | 2 | 0.4871 | 2 | 0.1138 | 2 | 0.0259 | 2 | 0.0059 | 2 |
| 4 | 1 | 0.9783 | 1 | 0.2366 | 1 | 0.0541 | 1 | 0.0122 | 1 |
| 12.5858 | 2 | 2.2871 | 1 | 0.5609 | 1 | 0.1309 | 1 | 0.0298 | 1 |
| 15.4142 | 2 | 2.3296 | 2 | 0.5751 | 1 | 0.1340 | 1 | 0.0304 | 1 |
| | | 2.3820 | 1 | 0.5902 | 2 | 0.1375 | 2 | 0.0313 | 2 |
| | | 2.8687 | 1 | 0.7301 | 1 | 0.1697 | 1 | 0.0386 | 1 |
| | | 3.0247 | 2 | 0.7332 | 2 | 0.1711 | 2 | 0.0389 | 2 |
| | | 3.4427 | 1 | 0.8313 | 1 | 0.1962 | 1 | 0.0446 | 1 |
| | | 4.4197 | 1 | 1.0996 | 1 | 0.2598 | 1 | 0.0591 | 1 |
| | | 4.4464 | 2 | 1.1241 | 2 | 0.2638 | 2 | 0.0599 | 2 |
| | | 4.6180 | 1 | 1.1709 | 1 | 0.2719 | 1 | 0.0616 | 1 |
| | | 5.2598 | 2 | 1.3120 | 1 | 0.3097 | 1 | 0.0705 | 1 |
| | | 5.3703 | 1 | 1.3435 | 2 | 0.3116 | 2 | 0.0707 | 2 |
| | | 5.4460 | 2 | 1.3576 | 1 | 0.3172 | 1 | 0.0720 | 1 |
| | | 5.4463 | 1 | 1.3945 | 2 | 0.3294 | 2 | 0.0751 | 2 |
| | | 5.4601 | 1 | 1.5662 | 1 | 0.3731 | 1 | 0.0853 | 1 |
| | | 6.6192 | 2 | 1.6627 | 2 | 0.3929 | 2 | 0.0895 | 2 |
| | | 6.7007 | 2 | 1.7328 | 2 | 0.4049 | 2 | 0.0920 | 2 |
| | | 7.8468 | 1 | 1.9195 | 1 | 0.4505 | 1 | 0.1026 | 1 |
| | | 7.8481 | 1 | 2.0174 | 1 | 0.4865 | 1 | 0.1114 | 1 |
| | | 8.4801 | 2 | 2.1509 | 2 | 0.5150 | 2 | 0.1179 | 2 |
| | | 8.5001 | 2 | 2.2085 | 2 | 0.5213 | 2 | 0.1187 | 2 |
| | | 8.7486 | 1 | 2.2149 | 1 | 0.5284 | 1 | 0.1204 | 1 |
| | | 8.8864 | 2 | 2.2517 | 1 | 0.5367 | 1 | 0.1231 | 1 |
| | | 9 | 1 | 2.2995 | 2 | 0.5469 | 2 | 0.1248 | 2 |
| | | 9.1163 | 2 | 2.3466 | 2 | 0.5675 | 2 | 0.1296 | 2 |
| | | 9.1668 | 2 | 2.3763 | 1 | 0.5683 | 2 | 0.1298 | 2 |
| | | 9.2832 | 1 | 2.3864 | 2 | 0.5712 | 2 | 0.1311 | 2 |
| | | 9.3041 | 2 | 2.3900 | 2 | 0.5806 | 1 | 0.1353 | 1 |
| | | 10.5892 | 2 | 2.3992 | 2 | 0.5897 | 2 | 0.1364 | 1 |
| | | 10.8516 | 2 | 2.4039 | 1 | 0.5905 | 1 | 0.1373 | 2 |
| | | 11.2427 | 2 | 2.4135 | 1 | 0.5915 | 1 | 0.1375 | 1 |
| | | 12.0270 | 2 | 2.4259 | 1 | 0.5990 | 1 | 0.1395 | 1 |
| | | 12.5055 | 2 | 2.4369 | 1 | 0.6014 | 2 | 0.1399 | 2 |
| | | 12.5145 | 2 | 2.4507 | 1 | 0.6145 | 1 | 0.1430 | 1 |
| | | 14.1299 | 2 | 2.4604 | 2 | 0.6145 | 1 | 0.1430 | 1 |
| | | 14.1407 | 2 | 2.4791 | 2 | 0.6183 | 1 | 0.1433 | 1 |
| | | 15.0161 | 2 | 2.5360 | 1 | 0.6248 | 2 | 0.1450 | 2 |
| | | 15.2155 | 2 | 2.6906 | 1 | 0.6668 | 1 | 0.1551 | 1 |
| | | | | 2.7170 | 2 | 0.6689 | 2 | 0.1558 | 2 |
| | | | | 2.7996 | 2 | 0.6815 | 2 | 0.1566 | 2 |
| | | | | 2.8159 | 1 | 0.6996 | 2 | 0.1607 | 2 |
| | | | | 2.8872 | 2 | 0.7034 | 1 | 0.1635 | 1 |

In order to guess the renormalization factor, we want to compare the ratio between the eigenvalues in different levels. It is difficult in general to decide if eigenvalues on different lines are associated with related eigenfunctions. Nevertheless, since the eigenfunctions corresponding to the eigenvalue of multiplicity 1 have good symmetry, as stated below, we can compare them easily.

**Proposition 5.1.** *Given an eigenfunction $u$ and the corresponding eigenvalue with multiplicity $1$, one the the following is true about $u$: invariant under $D_8$ (O), symmetric about a line through the mid-points of the opposite edges of an octagon (E), invariant under a $\frac{\pi}{4}$ rotation (R), or symmetric about a line through two opposite corners of an octagon (C).*

*Proof.* Let $E_\lambda$ be an $\lambda$-eigenspace of dimension 1 spanned by $u$. Then, by the construction of the Laplacian, one can easily see that for all $g \in D_8$, $\Delta(u \circ g) = \lambda u \circ g$. Since $\lambda$ has multiplicity 1, we have $u \circ g = \rho_g u$ for some $\rho_g \in \mathbb{C}$. Note that $u = u \circ g^8 = \rho_g^8 u$, and $u$ is not zero everywhere since it spans $E_\lambda$. Hence, $\rho_g^8 = 1$. In particular, $\rho_g \neq 0$. Since the matrix used to compute eigenfunction is symmetric, by spectral theorem, we can choose $u$ to be real-valued. Now pick $x \in POG$ such that $u(x) \neq 0$, we have $u \circ g(x) = \rho_g u(x) \neq 0$, which forces $\rho_g \in \mathbb{R}$ since $u$ is real-valued. Therefore $\rho_g = +1$ or $-1$. As $D_8$ is generated by reflection along the opposite edges of an octagon $e$ and rotation $r$, the symmetry of $u$ is completely determined by $\rho_e$ and $\rho_r$. One can easily show that the four possible cases $\rho_e = \rho_r = 1$, $\rho_e = 1$ and $\rho_r = -1$, $\rho_e = -1$ and $\rho_r = 1$, $\rho_e = \rho_r = -1$ correspond to the above cases respectively. □

*Remark* 5.2. Let $e, r$ and $c \in D_8$ be the reflection about a line through the mid-points of the opposite edges of an octagon, rotation anticlockwise by $\frac{\pi}{4}$ and reflection about a line through two opposite corners of an octagon respectively. If $u$ is of type (E), then $u \circ r = u \circ c = -u$. If $u$ is of type (R), then $u \circ e = u \circ c = -u$. If $u$ is of type (C), then $u \circ r = u \circ e = -u$.

*Remark* 5.3. There is an interesting observation related to the number of eigenvalues with multiplicity 1. Let $O_m, E_m, R_m$ and $C_m$ be the number of eigenvalues of type $O, E, R$ and $C$ respectively. We find that $O_m = C_m = 4 \times 8^{m-2} + 2^{m-2}$ and $E_m = R_m = 4 \times 8^{m-2} - 2^{m-2}$ for $m = 1, 2, 3$ and $4$.

Table 5.1 below shows a part of the comparison between the eigenvalues of multiplicity 1 corresponding to $-\triangle_1, -\triangle_2, -\triangle_3, -\triangle_4$ and $-\triangle_5$, which gives us the renormalization factor $r \approx 4.4$ when we consider the ratio of $R_4/R_5$. We therefore define the renormalized eigenvalues $\{r^m \lambda_k^m\}$ for $r = 4.4$. In Table 5.2 we show the renormalized eigenvalues on levels 3,4 and 5 in the beginning portion of the spectrum where the multiplicities are the same on these levels.

**Table 5.1 Symmetry Types and Ratios**

| Level1 | Level2 | Level3 | Level4 | Level5 | Type | $R_1/R_2$ | $R_2/R_3$ | $R_3/R_4$ | $R_4/R_5$ |
|---|---|---|---|---|---|---|---|---|---|
| 0 | 0 | 0 | 0 | 0 | O | | | | |
| 4 | 0.9783 | 0.2366 | 0.0541 | 0.0122 | C | 4.0887 | 4.1348 | 4.3734 | 4.4242 |
| | 2.2871 | 0.5609 | 0.1309 | 0.0298 | O | | 4.0776 | 4.2850 | 4.4000 |
| | 2.3820 | 0.5751 | 0.1340 | 0.0304 | E | | 4.1419 | 4.2918 | 4.4011 |
| | 2.8687 | 0.7301 | 0.1697 | 0.0386 | R | | 3.9292 | 4.3023 | 4.4013 |
| | 3.4427 | 0.8313 | 0.1962 | 0.0446 | C | | 4.1413 | 4.2370 | 4.3958 |
| | 4.4197 | 1.0996 | 0.2598 | 0.0591 | O | | 4.0194 | 4.2325 | 4.3956 |
| | 4.6180 | 1.1709 | 0.2719 | 0.0616 | E | | 3.9440 | 4.3064 | 4.4119 |
| | 5.3703 | | | | C | | | | |
| | 5.4463 | 1.3120 | 0.3097 | 0.0705 | O | | 4.1511 | 4.2364 | 4.3929 |
| | 5.4601 | 1.3576 | 0.3172 | 0.0720 | C | | 4.0219 | 4.2799 | 4.4036 |
| | | 1.5662 | 0.3731 | 0.0853 | C | | | 4.1978 | 4.3726 |
| | 7.8468 | 1.9195 | 0.4505 | 0.1026 | O | | 4.0879 | 4.2608 | 4.3926 |
| | 7.8481 | 2.0174 | 0.4865 | 0.1114 | R | | 3.8902 | 4.1468 | 4.3691 |

**Table 5.2 Renormalized Eigenvalues**

| 1 | 1 | 0 | 0 | 0 |
|---|---|---|---|---|
| 2,3 | 2 | 2.203168 | 2.2062656 | 2.21137664 |
| 4 | 1 | 4.580576 | 4.6084544 | 4.57267712 |
| 5 | 1 | 10.859024 | 11.1505856 | 11.16932608 |
| 6 | 1 | 11.133936 | 11.414656 | 11.39421184 |
| 7,8 | 2 | 11.426272 | 11.7128 | 11.73154048 |
| 9 | 1 | 14.134736 | 14.4557248 | 14.46765056 |
| 10,11 | 2 | 14.194752 | 14.5749824 | 14.58009344 |
| 12 | 1 | 16.093968 | 16.7131008 | 16.71650816 |
| 13 | 1 | 21.230176 | 22.1308032 | 22.15124736 |
| 14,15 | 2 | 21.762576 | 22.4715392 | 22.45109504 |
| 16 | 1 | 22.668624 | 23.1615296 | 23.08827136 |
| 17 | 1 | 25.40032 | 26.3814848 | 26.4240768 |
| 18,19 | 2 | 26.01016 | 26.9692544 | 26.49903872 |
| 20 | 1 | 26.283136 | 27.0203648 | 26.9862912 |
| 21,22 | 2 | 26.99752 | 28.0596096 | 28.14820096 |
| 23 | 1 | 30.321632 | 31.7821504 | 31.97125888 |
| 24,25 | 2 | 32.189872 | 33.4687936 | 33.5454592 |
| 26,27 | 2 | 33.547008 | 34.4910016 | 34.4824832 |
| 28 | 1 | 37.16152 | 38.375392 | 38.45546496 |
| 29 | 1 | 39.056864 | 41.442016 | 41.75378944 |
| 30,31 | 2 | 41.641424 | 43.86976 | 44.19005184 |
| 32,33 | 2 | 42.75656 | 44.4064192 | 44.48989952 |
| 34 | 1 | 42.880464 | 45.0112256 | 45.12707584 |
| 35 | 1 | 48.72912 | 45.7182528 | 46.13906176 |
| 36,37 | 2 | 44.51832 | 46.5871296 | 46.77623808 |
| 38,39 | 2 | 45.430176 | 48.34192 | 48.57532416 |

## 6. Eigenfunctions of Multiplicity 1

**Figures 6.0 - 6.8 Multiplicity 1 Eigenfunctions**

Figure 6.0 4th Eigenfunction- C    $\lambda = 4, 4.1088, 4.07425, 3.8195$

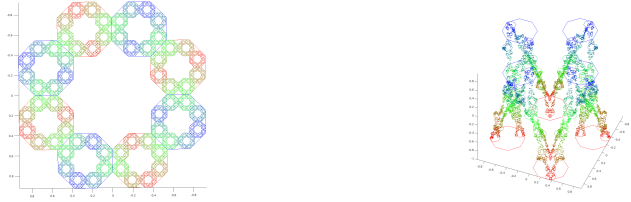

Figure 6.1 5th Eigenfunction- O    $\lambda = 9.60582, 9.65898, 9.2418$

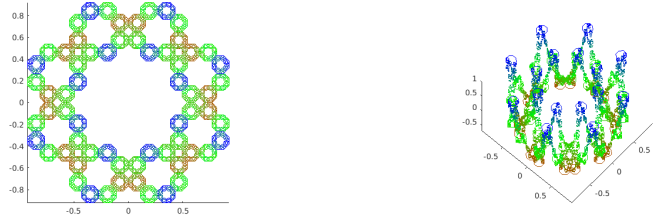

Figure 6.2 6th Eigenfunction- E    $\lambda = 10.0044, 9.903222, 9.460668$

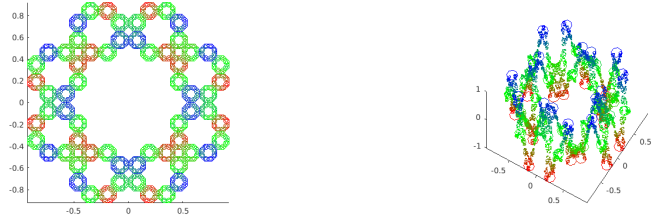

Figure 6.3 9th Eigenfunction- R    $\lambda = 12.04854, 12.5723, 11.9811$

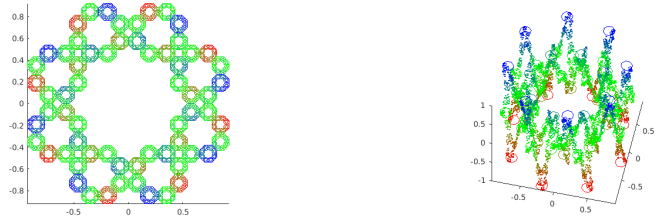

Figure 6.4 12th Eigenfunction- C    $\lambda = 14.45934, 14.31498, 13.852112$

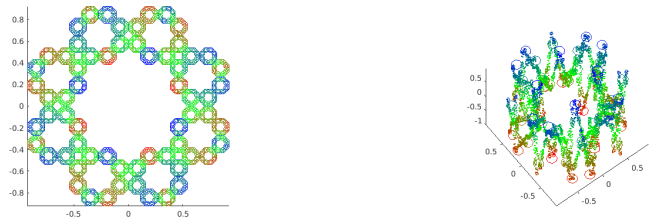

Figure 6.5 13th Eigenfunction $\lambda = 18.56274, 22.59264, 21.865439$

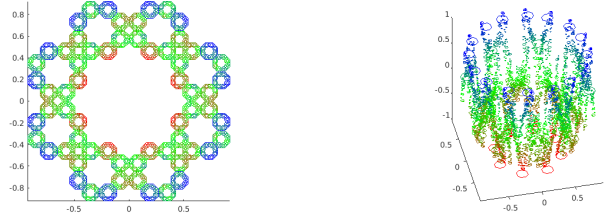

Figure 6.6 16th Eigenfunction $\lambda = 19.3956, 20.162898, 19.1966838$

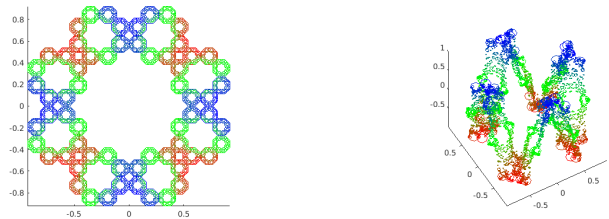

Figure 6.7 17th Eigenfunction $\lambda = 22.87446, 22.59264, 21.8654394$

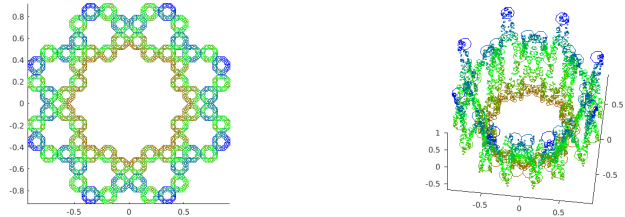

Figure 6.8 20th Eigenfunction $\lambda = 22.55526, 23.377872, 22.3949544$

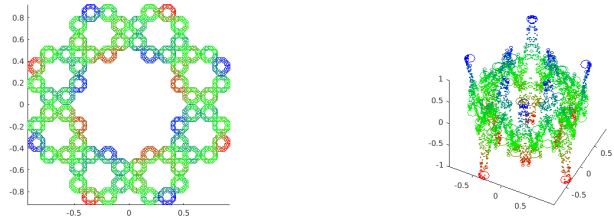

Figure 6.0-6.8 are the pictures for the eigenfunctions corresponding to some of the eigenvalues listed above. It has an overlay of the eigenfunctions lined up on consecutive levels. When four $\lambda$ values are given this indicates eigenvalues for levels 1,2,3,4. When there are only three values given this indicates eigenvalues for level 2,3,4. Blue indicates a positive value, red indicates a negative value and green indicates a value close to 0. Each eigenfunction is labelled $C, O, E$ or $R$ depending on the symmetry type.

Let $N_m(x)$ be the eigenvalue counting function for $\triangle_m$, defined to be the number of eigenvalues less than or equal to $x$. It is of interest to study the growth rate of $N(x)$ and the gaps between the eigenvalues. Figures 6.9 - 6.13 are the graphs of $x$ against $N(x)$ for level 1, 2 ,3 ,4 and 5 respectively.

Figure 6.9 Level 1

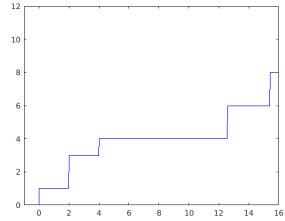

Figure 6.10 Level 2

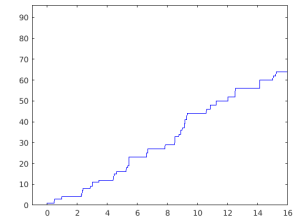

Figure 6.11 Level 3

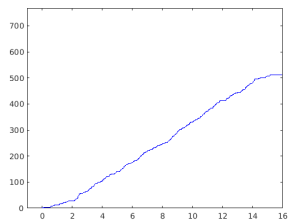

Figure 6.12 Level 4

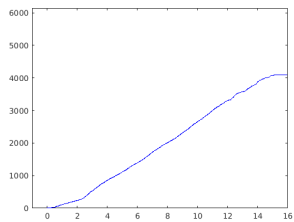

Figure 6.13 Level 5

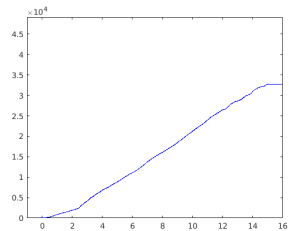

Next we examine the Weyl ratios for $\triangle_m$ defined to be $N_m(t)/t^{\{}\alpha$ where $\alpha$ is chosen to give the best fit. In Figures 6.14-6.18 we display the log-log plots for $m = 1, 2, 3, 4, 5$ respectively.

Figure 6.14 Level 1: $\alpha = 0.4863$

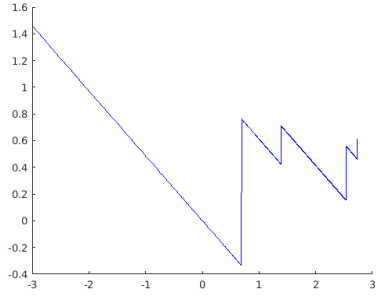

Figure 6.15 Level 2: $\alpha = 1.0615$

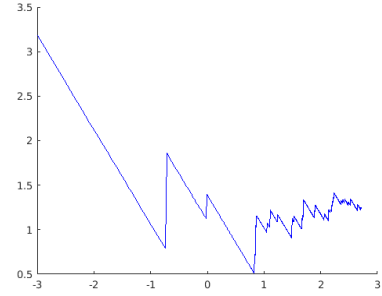

Figure 6.16 Level 3: $\alpha = 1.3181$

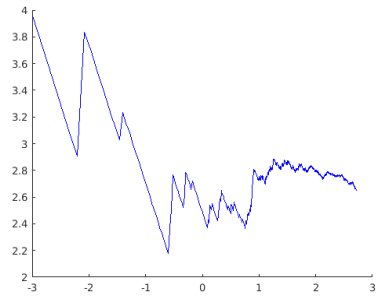

Figure 6.17 Level 4: $\alpha = 1.3957$

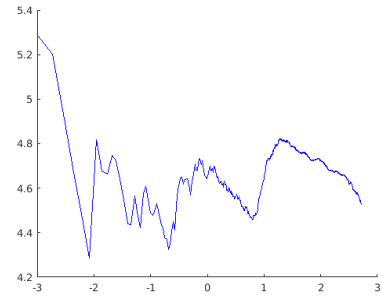

Figure 6.18 Level 5: $\alpha = 1.4105$

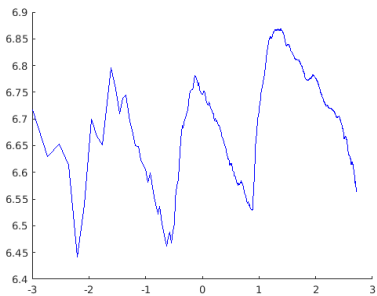

The $\alpha$'s for level 1, 2, 3, 4, 5 are 0.4863, 1.0615, 1.3181, 1.3957, 1.4105 respectively. Therefore, we would expect that the $N(x) \geq O(x^{1.4})$ for the eigenvalue counting function in the projective octagasket..

## 7. Heat Equation

The heat equation

$$\frac{\partial(u(x,t))}{\partial t} = \triangle u(x,t)$$

on our graph approximations given initial conditions

$$u(x,0) = f(x)$$

can be solved by using the heat kernel, a construction which utilizes an orthonormal basis of eigenfunctions for the Laplacian. The eigenfunctions we found at each level can be renormalized to construct such a basis.

Given an orthonormal eigenbasis $\{\phi_i\}$ the heat kernel is defined as

$$h_t(x,y) = \sum_\lambda e^{-\lambda t} P_\lambda(x,y)$$

where $P_\lambda(x,y) = \sum_i \phi_i(x)\phi_i(y)$ Now, given our initial heat distribution $f(x)$ we can solve for $u(x,t)$ for any time $t$ by simply taking the integral of the product of the heat kernel and the initial $f(x)$ with respect to our measure. The solution to the heat equation is given as

$$u(x,t) = \int_K h_t(x,y) f(y) d\mu(y)$$

Of course since we only have approximations to work with our integral turns into a finite sum. We were able to find solutions to heat equation especially when setting $f(x_0) = 1$ and $f(x) = 0$ for all $x \neq x_0$. This allows us to simplify the finite sum to $u(x,t) = h_t(x, x_0)$.

Figure 7.0 - 7.5 show the solution to the heat equation when setting $f(00...0) = 1$ and $f(x) = 0$ if $x \neq 00...0$.

The Heat Kernel at varying time $t$

Figure 7.0: t= 0       Figure 7.1: t= 0.1       Figure 7.2: t = 0.2

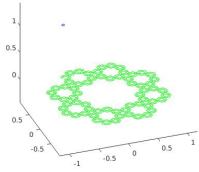 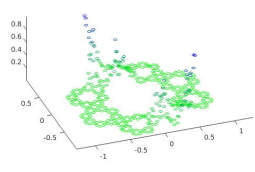 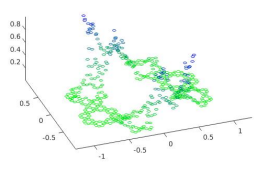

Figure 7.3: t= 0.3       Figure 7.4: t= 0.4       Figure 7.5: t = 1

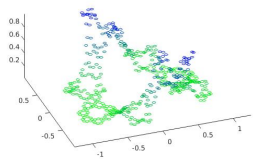 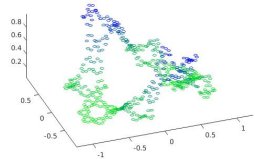 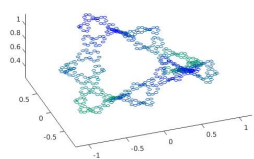

## 8. Wave Equation

Similar to the heat equation, we can solve the wave equation

$$\frac{\partial^2 u(x,t)}{\partial t^2} + \triangle u(x,t) = 0$$

on our graph approximations given the initial conditions

$$\begin{cases} u(x,0) = f_0(x) \\ \frac{\partial u}{\partial t}(x,0) = f_1(x) \end{cases}$$

If $f_0 = 0$, then the solution is given by

$$u(x,t) = \sum_\lambda \frac{\sin t\sqrt{\lambda}}{\sqrt{\lambda}} \int_K P_\lambda(x,y) f_1(y) d\mu(y)$$

, where $P_\lambda(x,y)$ is defined as the same as before.

We us use a level-3 cell graph approximation of the POG as an illustration for the solution to the wave equation. If we assume $f_1(x) = 0$ everywhere, except $f_1(000) = 1$ and $f_1(001) = -1$ then we obtain the solutions given in Figures 8.0-8.19

**Solution to the Wave Equation for varying time $t$**

Figure 8.0: t = 0 — Figure 8.1: t = 1

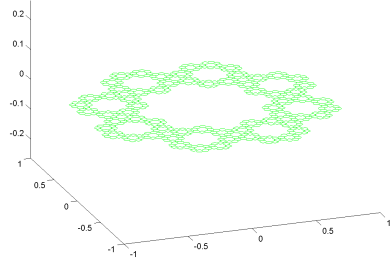
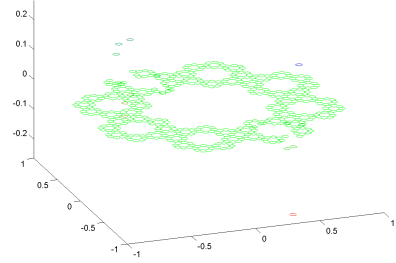

Figure 8.2: t = 2 — Figure 8.3: t = 3

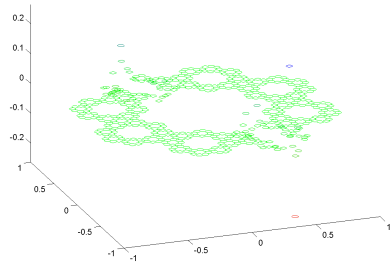
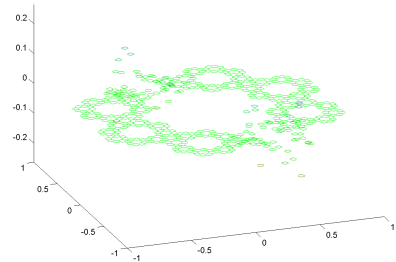

Figure 8.4: t = 4 — Figure 8.5: t = 5

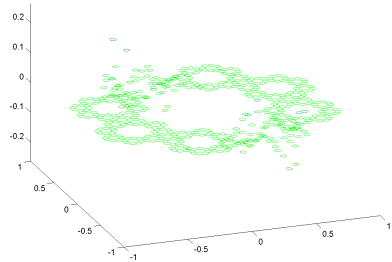
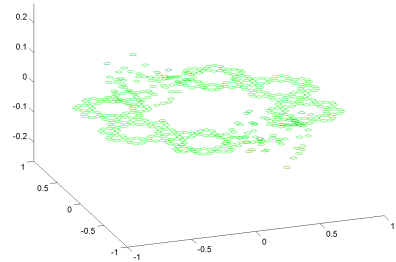

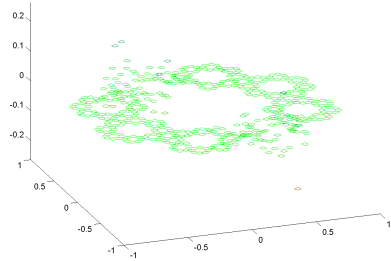

Figure 8.6: t = 6

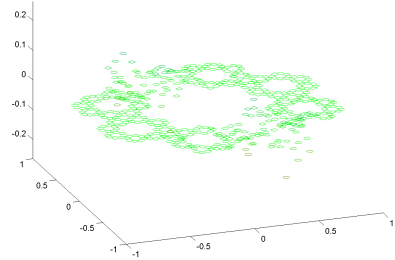

Figure 8.7: t = 7

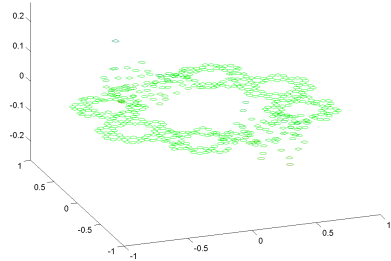

Figure 8.8: t = 8

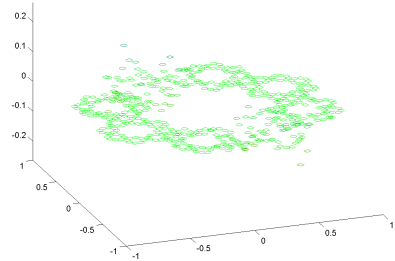

Figure 8.9: t = 9

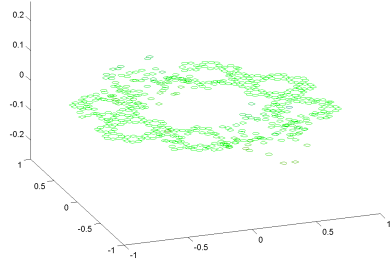

Figure 8.10: t = 10

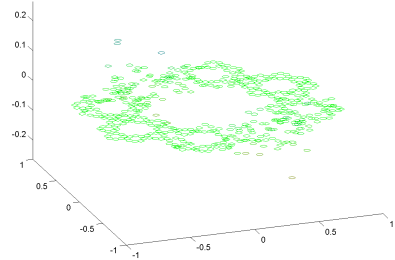

Figure 8.11: t = 11

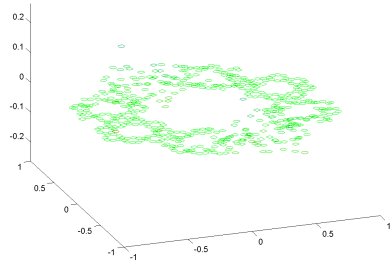

Figure 8.12: t = 12

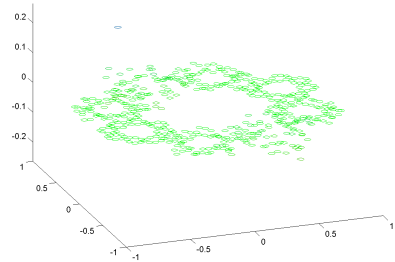

Figure 8.13: t = 13

Figure 8.14: t = 14     Figure 8.15: t = 15

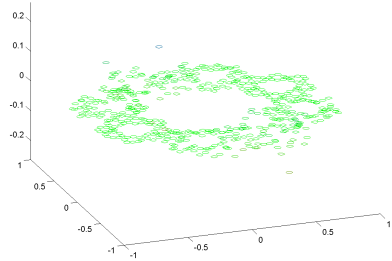 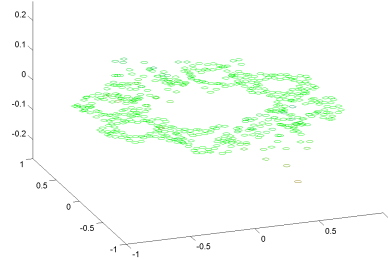

Figure 8.16: t = 16     Figure 8.17: t = 17

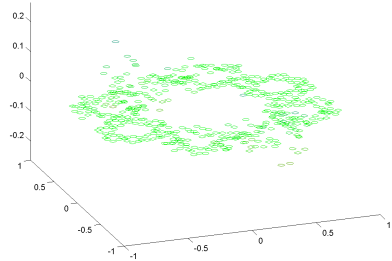 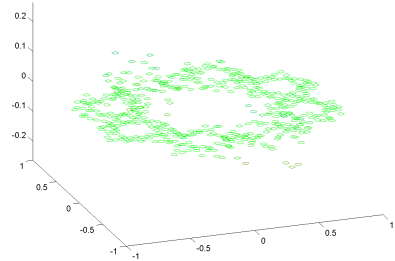

Figure 8.18: t = 18     Figure 8.19: t = 19

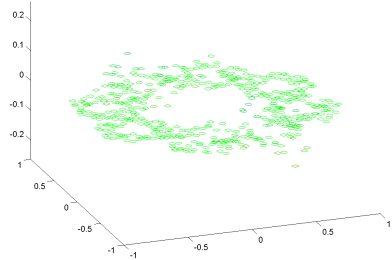 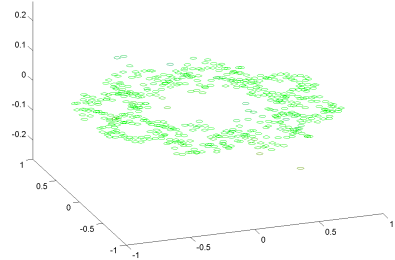

## 9. Geometry

### 9.1. Metric in the Finite Graph.

Given $\Gamma_m$, there exists a natural metric $D_m$ on it. We call $l$ a path in $\Gamma_m$ if $l$ is a sequence of $m$-cells $c_0 c_1 c_2 ... c_n$ with each pair of consecutive cells $(c_i, c_{i+1})$ is adjacent but not equal. We call $n$ the length of $l$. Let $w, v$ be $m$-cells in $\Gamma_m$. Then we define $D_m(w,v)$ to be the minimum of the length of all paths with end points $w$ and $v$. Let $d_m = diam(\Gamma_m)$. Our experimental results show that $d_1 = 2, d_2 = 7, d_3 = 15$, $d_4 = 29$ and $d_5 = 57$, and $d_m$ is attained if $w = X00...0$, $v = (X \pm 2)00...0$ for some $X \in \mathbb{Z}_8$, which gives us the following conjecture:

**Conjecture 9.1.** $d_m = 2^{m-2} \times 7 + 1 \; for \; m > 2$

Though we cannot prove it at this stage, we have the following propositions instead.

**Proposition 9.2.** $d_{m+1} \leq 2d_m + 2^m + 1$ for all $m > 0$

*Proof.* Let $u = u_1u_2...u_{m+1}, v = v_1v_2...v_{m+1}$ be $(m+1)$-cells in $\Gamma_m$. Let $U$ and $V$ be the 1-cell containing $u$ and $v$ respectively.

If $U = V$ (i.e $u_1 = v_1$), then obviously, by the definition of $d_m$, $D_{m+1}(u,v) \leq d_m$.

If the 1-cell containing $u$ and the 1-cell containing $v$ are adjacent (i.e. $u_1 = v_1 \pm 1$ or $v_1 \pm 4$), one can easily construct a path $l$ starting from $u$, passing through the boundary between $U$ and $V$, and then ending at $v$ with $length(l) \leq 2d_m + 1$. So, $D_{m+1}(u,v) \leq 2d_m + 1$

If $u_1 = v_1 \pm 2$, without loss of generality, say $u_1 = v_1 - 2$, then one can construct a path $l_1$ from $u$ to $u_1 200...0$, a path $l_2$ from $(u_1+1)600...0$ to $(u_1+1)200...0$, and a path $l_3$ from $(u_1+2)600...0$ to $v$ with $length(l_1), length(l_3) \leq d_m$ and $length(l_2) \leq 2^m - 1$ (in fact it suffices to take $l_2$ to be the straight line from $(u_1+1)600...0$ to $(u_1+1)200...0$ in the Euclidean sense). Connecting them, we obtain a path from $u$ to $v$ with length $\leq 2d_m + 2^m + 1$. Hence, $D_{m+1}(u,v) \leq 2d_m + 2^m + 1$.

If $u_1 = v_1 \pm 3$, without loss of generality, say $u_1 = v_1 - 3$, then one can construct a path $l_1$ from $u$ to $u_1 600...0$, and a path $l_2$ from $(u_1 - 3)600...0$ with $length(l_1), length(l_2) \leq d_m$. Note that $(u_1+4)6000...0$ is adjacent to both $u_1 600...0$ and $(u_1+3)200...0$. Connecting them, we obtain a path from $u$ to $v$ with length $\leq 2d_m + 2$. Hence, $D_{m+1}(u,v) \leq 2d_m + 2$

Figures 9.0-9.3 show how one can construct a path from $u_1u_2...u_{m+1}$ to $u_1x_2...x_{m+1}$, $(u_1+1)x_2...x_{m+1}$, $(u_1+2)x_2...x_{m+1}$ and $(u_1+3)x_2...x_{m+1}$ respectively. The gradient from red to blue indicates the distance from the starting 1 cell. Since the maximal distance covered in each 1 cell is bounded by $d_m$ we can get a better bound on $d_{m+1}$.

It follows that $d_{m+1} \leq 2d_m + 2^m + 1$. □

Figure 9.0 Path 1

Figure 9.1 Path 2

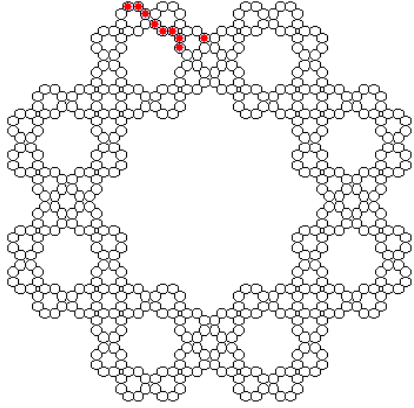
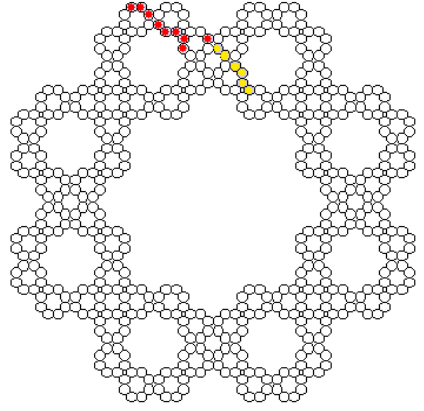

Figure 9.2 Path 3

Figure 9.3 Path 4

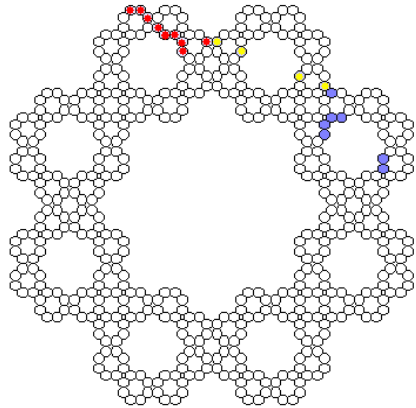
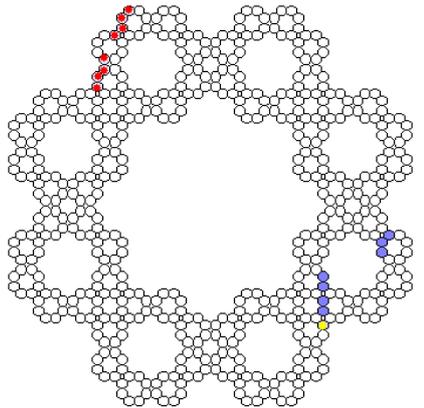

**Corollary 9.3.** $d_m \leq 2^{m-1}(m+2) - 1$ for all $m > 0$

*Proof.* Directly follows from Proposition 9.2 and the fact that $d_1 = 2$. □

**Proposition 9.4.** $d_m \geq 2^{m-1} - 1$ for all $m > 0$

*Proof.* It is trivial for $m = 1$.

For $m > 1$, first observe that $d_m \geq D_m(\text{outer boundary}, \text{inner boundary})$. What we need is in fact $D_m(\text{outer boundary}, \text{inner boundary}) \geq 2^{m-1} - 1$. We will construct layers of cells $\{\mathcal{C}_{m,k}\}_{k=1}^{2^{m-1}}$ in $\Gamma_m$ such that $\mathcal{C}_{m,1}$ is the inner boundary, $\mathcal{C}_{m,2^{m-1}}$ is the outer boundary. $\cup_{U \in \mathcal{C}_{m,i}} U$ only intersect $\cup_{U \in \mathcal{C}_{m,i+1}} U$ and $\cup_{U \in \mathcal{C}_{m,i-1}} U$ (if exists).

For m=2, we show in Figure 9.4 that it suffices to take the following layers, grey denotes $\mathcal{C}_{2,1}$ and black denotes $\mathcal{C}_{2,2}$.

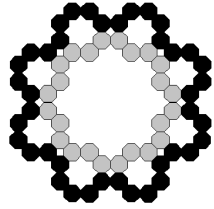

**Figure 9.4 Outer Layer: Level 2**

Suppose the statement is true for $m = n$. We shall split $\mathcal{C}_{m,i}$ into two layers in $\Gamma_{n+1}$. $\mathcal{C}_{n+1,1}$ and $\mathcal{C}_{n+1,2^n}$ can be obtained by modifying the identification algorithm in section 4. $\mathcal{C}_{n+1,2}$ and $\mathcal{C}_{n+1,2^n-1}$ will be the collection of the remaining $(n+1)$-cells with the same prefix as the cells in $\mathcal{C}_{n,1}$ and $\mathcal{C}_{n,2^{n-1}})$ respectively.

To construct layer $\mathcal{C}_{n+1,2i-1}$, one first draw a piecewise straight line $l$ to connect the centres of every pair of adjacent cells in $\mathcal{C}_{n,i}$. Now replace the octagons in $\mathcal{C}_{n,i}$ by $\Gamma_1$. One can inductively define $\mathcal{C}_{n+1,2i-1}$ by including a cell adjacent to some cells in $\mathcal{C}_{n+1,2i-2}$ and the cell clockwise adjacent to the previous cell, unless the previous cell touches the line $l$. If this is the case, one will move to the next copy of $\Gamma_1$ and the cell adjacent (in the sense of the octagasket) to the previous cell and repeat the above procedures, until the cells start to repeat.

Then, $\mathcal{C}_{n+1,2i}$ will be the collection of the remaining $(m+n+1)$-cells with the same prefix as the cells in $\mathcal{C}_{n,i}$. Figures 9.5-9.11 show how to construct the layers from $m = 2$ to $m = 3$.

Figure 9.5 Layer 0

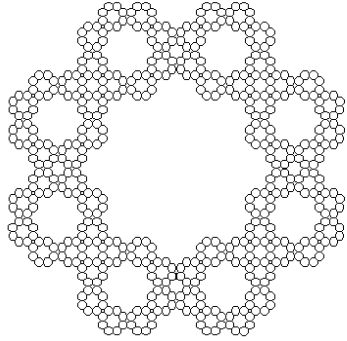

Figure 9.6 Layer 1

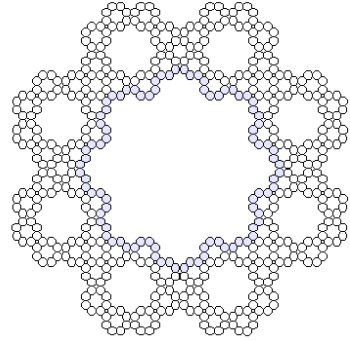

Figure 9.7 Layer 2

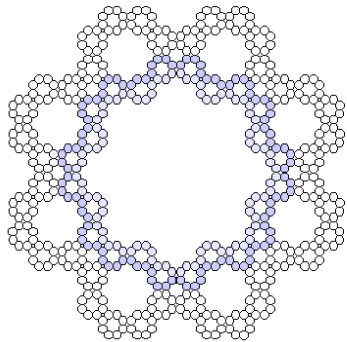

Figure 9.8 From layer 2 to layer 3

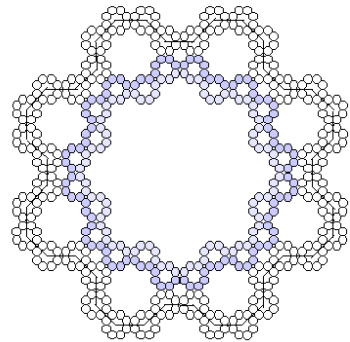

Figure 9.9 From layer 2 to layer 3

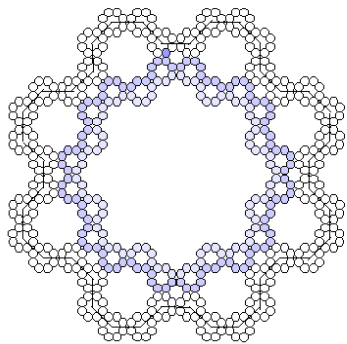

Figure 9.10 Layer 3

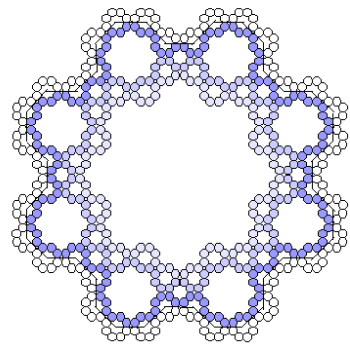

Figure 9.11 Layer 4

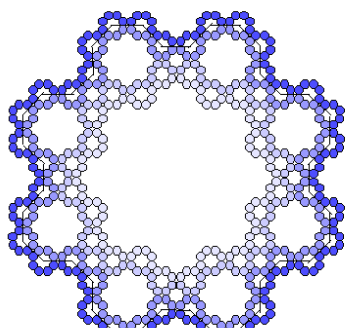

One can observe that the layers constructed satisfy our requirement. Now given any path from the inner boundary to the outer boundary, it must pass through each layer, which requires the path of length at least $2^n - 1$, so our statement follows. □

Figures 9.12 - 9.15 display the distance between 00...0 and other $m$-cells in $\Gamma_m$. The darkness of the cell represents the distance between that cell and 00...0, and the red cell is 00...0. The similarity of these pictures gives us evidence for the existence of the metric $D : POG \times POG \to [0, 1]$, where $D$ is determined by

$$D(v,w) = \lim_{m \to \infty} \frac{D_m(V_m, W_m)}{d_m}$$

for all $v, w \in V_*$, where $V_m$ and $W_m$ are $m$-cells in $\Gamma_m$ with $\bigcap_{n=1}^{\infty} V_n = \{v\}$ and $\bigcap_{n=1}^{\infty} W_n = \{w\}$. In addition, it is believed that the metric induces the same topology as the usual quotient topology.

Figure 9.12 Distance $\Gamma_1$

Figure 9.13 Distance $\Gamma_2$

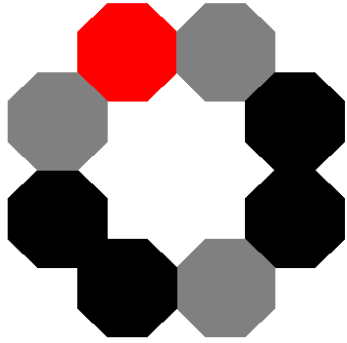

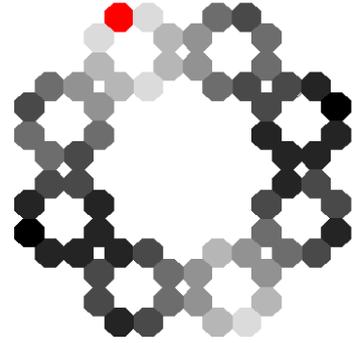

Figure 9.14 Distance $\Gamma_3$

Figure 9.15 Distance $\Gamma_4$

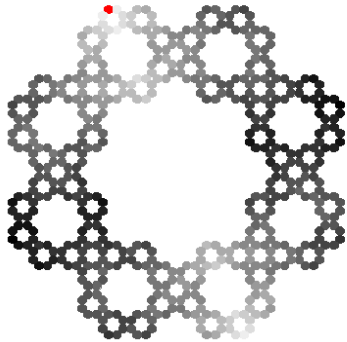

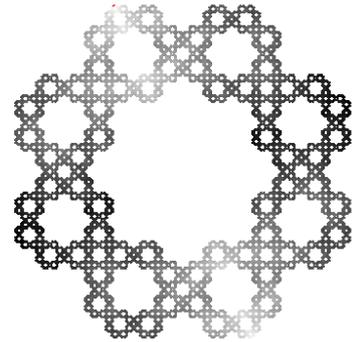

### 9.2. Cardinality of Metric Balls in Finite Graph.

The existence of metric allows us to discuss balls of radius $n$ of an $m$-cell $v$ in $\Gamma_m$. The experimental results show that the cardinality of $B_m(v, n)$ for $m = 2, 3, 4, 5$ are $C_2 n^{2.5377}, C_3 n^{2.7967}, C_4 n^{2.9483}, C_5 n^{3.0203}$ approximately for some constants $C_2, C_3, C_4, C_5$ when $r$ is not too large. In fact, using the similar technique developed by [?], we have the following bound: $\#B_{m+1}(v,n) \leq \sqrt[3]{32}n^3$ if $n < d_m$.

Before getting into the bound, we first give some definitions and develop some lemmas. We first define an $m$-edge in $\Gamma_{m+k}$. Basically, an $m$-edge is a path travelling from a corner point of an octagon to the adjacent corner point with its locus is totally contained in the outer boundary of the $m$-cell, where the octagon is the $m$-cell in $\Gamma_m$ with the same address as that $m$-cell in $\Gamma_{m+k}$.

**Definition 9.5.** In $\Gamma_{m+k}$ with $m > 0$, an $m$-edges $E$ is a piecewise straight line in the outer boundary of the $m$-cell $X_1X_2...X_m$ travelling from the point $X_1X_2...X_m00...0$ to $X_1X_2...(X_m+1)00...0$. We write $E \in E_{vw}$ if $E \subseteq v \cap w$. If $E = v \cap w$, we write $E = E_{vw}$. A cell $u$ is on an $m$-edge $E$ if $u \cap E \neq \emptyset$. Two edges $E, E'$ are adjacent if $E \cap E' \neq \emptyset$.

*Remark* 9.6. Each $m$-cell is on 8 distinct $m$-edge. Two $m$-cells may intersect on more than one $m$-edges.

**Lemma 9.7.** *In $\Gamma_{m+k}$ with $m, k > 0$, if a path, with length $< 2^k - 1$, lies in an $m$-cell and begins at an $(m+k)$-cell on the $m$-edge, then it cannot pass through nor end at any $(m+k)$-cell on the non-adjacent $m$-edge.*

*Proof.* The proof is similar to Proposition 8.4, which is basically by induction on $k$. Without loss of generality, we first assume the $m$-edge is from the top left corner to the top right corner. Let $\mathcal{C}_{k,1}$ be the collection of the $(m+k)$-cell on that $m$-edge. We will construct layers of cells $\{\mathcal{C}_{k,i}\}_{i=1}^{2^k}$ in the $m$-cell such that the $\cup_{U \in C_{k,i}} U$ only intersect $\cup_{U \in C_{k,i+1}} U$ and $\cup_{U \in C_{k,i-1}} U$

For $k = 1$, it is suffices to take the layers in Figure 9.16, where grey denotes the cells in $\mathcal{C}_1$, black denotes the cells in $\mathcal{C}_2$.

**Figure 9.16 Cells in $\mathcal{C}_1$ and $\mathcal{C}_2$**

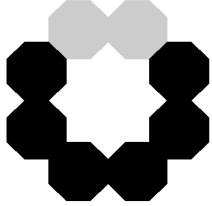

Suppose the statement is true for $k = n$. We split the $\mathcal{C}_{n,i}$ into two layers in $\Gamma_{m+n+1}$. $\mathcal{C}_{n+1,1}$ is clear by the identification algorithm in section 4. $\mathcal{C}_{n+1,2}$ will be the collection of the remaining $m+n+1$-cells with the same prefix as the cells in $\mathcal{C}_{n,1}$.

To construct layer $\mathcal{C}_{n+1,2i-1}$, one first draw a piecewise straight line $l$ to connect the centres of every pair of adjacent cells in $\mathcal{C}_{n,i}$. Now replace the octagons in $\mathcal{C}_{n,i}$ by $\Gamma_1$. One can inductively define $\mathcal{C}_{n+1,2i-1}$ by including the most top-left cell adjacent to some cells in $\mathcal{C}_{n+1,2i-2}$ and the cell clockwisely adjacent to the previous cell, unless the previous cell touches the line $l$. If this is the case, one will move to the next copy of $\Gamma_1$, and the cell adjacent (in the sense of octagasket) to the previous cell and repeat the above procedures, until reach the most top-right cell which is adjacent to some cells in $\mathcal{C}_{n+1,2i-2}$.

Then, $\mathcal{C}_{n+1,2i}$ will be the collection of the remaining $(m+n+1)$-cells with the same prefix as the cells in $\mathcal{C}_{n,i}$. Figures 9.17 - 9.22 show how to construct the layers from $k=1$ to $k=2$.

Figure 9.17 Layer 0     Figure 9.18 Layer 1     Figure 9.19 Layer 2

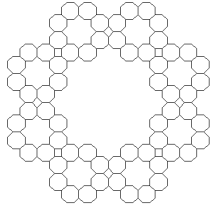 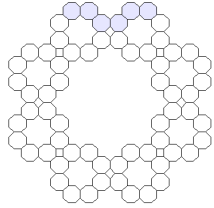 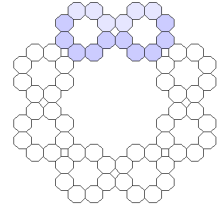

Figure 9.20 Intermediate Step     Figure 9.21 Layer 3     Figure 9.22 Layer 4

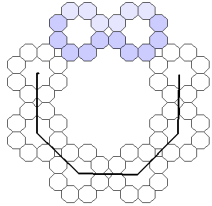 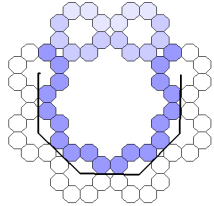 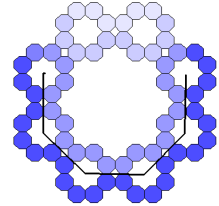

One can see that the layers constructed satisfy our requirement, only the cells on $C_{n+1,2^{n+1}}$ are on the non-adjacent $m$-edge. Now given any path from the original $m$-edge to any non-adjacent $m$-edge, it must across every layer if it does not go outside the $m$-cell, which require the path of at least length $2^k - 1$, so our statement follows. □

**Definition 9.8.** Consider $\Gamma_{m+k}$ with $m, k > 0$. Write a path $l$ in $\Gamma_{m+k}$ as a sequence of $m+k$-cells $c_0, c_1, ..., c_n$, with each $c_i$ has address $c_{i,1}c_{i,2}...c_{i,m+k}$. We can form a sequence of $m$-cells $\{V_i\}$ associate with $l$ as follows:
  (1) For each $c_i$, let $V_i = c_{i,1}c_{i,2}...c_{i,m}$.
  (2) While consecutive elements in the new sequence are equal, delete all but one of them, and relabel the index if necessary.

We call the set $\{V_i\}$ the $m$-sequence of $l$.

**Definition 9.9.** Given an $m$-sequence $\{V_i\}$ for a path $l$ in $\Gamma_{m+k}$, we call a tuple $(i, X, Y)$ an $m$-segment-of-two if $X = V_i, Y = V_{i+1}$ and $V_{i-1} \neq X$ and $Y$ (if it exists).

A path $l$ is said to enter an $m$-segment-of-two $(i, X, Y)$ through an $m$-edge $E$ if there exists some consecutive cells $c_n, c_{[}n+1$ such that both of them are on $E$ and only $c_{n+1}$ is within $(i, X, Y)$. The notion of a path exiting is analogous, with $c_n$ in $(i, X, Y)$ instead. The collection of the $m$-edges along which the path enters is denoted by $Ent(i)$, and the collection of the $m$-edges along which it exists is denoted by $Exit(i)$.

**Lemma 9.10.** Let $(i, X, Y)$ be an $m$-segment-of-two associated to a path $l$ of length at most $2^k - 1$.

(1) Let $E_1 \in Ent(i)$ (if exists). $E_1 \notin E_{XY}$, but there exists some $E_2 \in E_{XY}$ such that $E_1$ intersects $E_2$. Similar results also hold for $Exit(i)$.
(2) $Ent(i)$ and $Exit(i)$ are singleton (if exist).
(3) The edge that $l$ passes through $E_{XY}$ is uniquely determined by $Ent(i)$.
(4) If $Ent(i)$ and $Exit(i)$ exists, then there exists unique $E_1 \in Ent(i), E_2 \in E_{XY}, E_3 \in Exit(i)$ such that $E_1 \cap E_2 \cap E_3$ is a singleton, which is called the centre of the $(i, X, Y)$.

*Proof.* For the first statement, By the definitions of $m$-segments-of-two and $Ent(i)$, $E_1 \notin E_{XY}$. Let $E_2 \in E_{XY}$ be the $m$-edge such that $l$ passes through. Since the path at most $2^k - 1$ steps, and one step is required to enter the $m$-segment-of-twom, Lemma 8.7 forces $E_1$ is adjacent to $E_2$, i.e. $E_1$ intersects $E_2$. Similar results hold for $Exit(i)$.

Note that $Ent(i)$ is a singleton. If not, then there exists some $E' \in Ent(i)$ such that $E' \neq E_1$. By similar reasoning, we have $E'$ intersects $E_2$, in particular, $E_1$ and $E'$ are non-adjacent to each other, i.e. there is a cell on two non-adjacent edge simultaneously, which is impossible. Similarly, we know that $E_2$ is determined by $E_1$.

For the forth statement, by similar reasoning, there exists some adjacent $E_3 \in Exit_i$ such that $E_2$ intersects $E_3$. Therefore, there are two possibilities for $E_3$. One intersects $E_1$ and one is not, as shown in Figure 9.23.

**Figure 9.23 Possibilities for $E_i$**

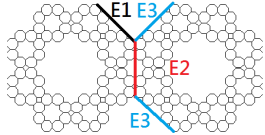

Similar to Lemma 8.7, we can construct layers to show this is the only possible way. Given 2 copies of $\Gamma_{k-1}$, we can identify one of the edge and induct on $k$ with the initial 2 copies of $\Gamma_1$ where $\mathcal{C}_{k,1}$ is taken to be the cells on $E_1$ and the $E_4$ as shown below in Figure 9.24 and Figure 9.25. Then, the length of $l$ forces $E_3 = E_4$. The existence of the centre follows.

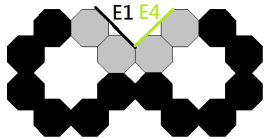

**Figure 9.24**

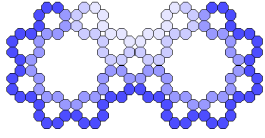

**Figure 9.25**

$\square$

**Lemma 9.11.** *If $l$ is a path of length at most $2^k - 1$ and at least two $m$-segments-of-two, then the centres of the $m$-segments-of-two are the same.*

*Proof.* First we prove the case when there is two $m$-segments-of-two. Let $(i, V, W)$ and $(j, X, Y)$ be two consecutive $m$-segments-of-two (i.e. $X \in \{V, W\}$ but $Y \notin \{V, W\}$), $E_{1,i}, E_{2,i}, E_{3,i}, E_{1,j}, E_{2,j}$ and $E_{3,j}$ be the $E_1, E_2$ and $E_3$ of $(i, V, W)$ and $(j, X, Y)$ respectively. It suffices to prove that $E_{2,i} = E_{1,j}$ and $E_{3,i} = E_{2,j}$. Note that $X = V$ or $W$ and $(i, V, W)$, $(j, X, Y)$ are consecutive $m$-segments-of-two. Hence, $E_{1,j} \in E_{VW}$. Also note that the edge that $l$ passes through $E_{VW}$ is unique, which forces $E_{1,j} = E_{2,i}$. Similarly, we have $E_{3,i} = E_{2,j}$. Since $E_{2,i} \cap E_{3,i} = E_{1,j} \cap E_{2,j}$ is in fact a singleton, thus the centres of the $m$-segments-of-two must be that singleton.

The remaining cases follows from induction. $\square$

**Proposition 9.12.** *For a path $l$ of length at most $2^k - 1$. Let $\{V_i\}_{i=1}^n$ be the $m$-sequence associated with $l$. Then there exists a point $x$ inside $V_i$ for all $i$.*

*Proof.* If $n = 1$, then it is trivial. Otherwise, it suffices to take the common centres to be $x$. $\square$

**Corollary 9.13.** *Let $x \in \Gamma_{m+k}$, $V_x$ be a $k$-cell containing $x$, $\mathcal{V}$ be the collection of all $k$-cells adjacent to $V_x$ with $V_x \in \mathcal{V}$. Any path $l$ starts at $V_x$ with length at most $2^m - 1$ cannot leave $\bigcup_{V \in \mathcal{V}} V$; i.e. $B_{m+k}(V_x, 2^m - 1) \subseteq \bigcup_{V \in \mathcal{V}} V$.*

*Proof.* Let $\{V_i\}_{i=1}^n$ be the $m$-sequence associated with $l$ with $V_x = V_1$. By proposition 9.12, every $V_i$ is adjacent to $V_x$. $\square$

**Lemma 9.14.** $\#\mathcal{V} \leq 25$

*Proof.* Each $k$-cell $v$ has exactly 8 corner points. Before identification, each corner point lies in at most 2 different cells at the same time. Hence, it will lies in at most 4 cells at the same time after identification, and at most 3 are different from $v$. Hence, it can be adjacent to at most 25 $k$-cells, 24 cells are other cells, and the left one cell is $v$. $\square$

**Proposition 9.15.** $\#B_m(v, n) \leq 200n^3$ *if* $n < 2^{m-1}$.

*Proof.* Let $l$ be a path starting at $v$ and with length $\leq n$. Let $k \in \mathbb{N}$ such that $2^{k-1} \leq n \leq 2^k - 1$. Then, $\#B_m(v, n) \leq \#B_m(v, 2^k - 1) \leq \#\mathcal{V} \times \#_{m-k} \leq 25 \times 8^k \leq 200n^3$, where $\#_k$ is the number of $m$-cells lying in a $k$-cell. $\square$

## 10. Discussion

In this paper we have computed the spectrum of the renormalized Laplacian on the level m approximations $POG_m$ for $m \leq 5$. We are interested in the Laplacian on $POG$ that is the limit of these approximations. We believe that this yields strong experimental evidence of the existence of the limiting Laplacian. Table 5.2 shows that with the choice of $r = 4.4$ we have a fairly accurate approximation to the limiting eigenvalues, but we do not have a clear conjecture for the rate of convergence. Figure 6.0-6.8 give reasonable evidence that the associated eigenfunctions converge. The fluctuations in the approximations are of the same order of magnitude as are seen in other fractal Laplacians where the existence of a limit has been proven. In Figure 6.9-6.13 we have graphs of the eigenvalue counting function on levels 1-5. This sequence of graphs seems to be converging, with no visible difference between levels 4 and 5. We also have graphs of the Weyl ratios $N(t)/t^\alpha$ for the best choice of $\alpha$ for each level. Here we do not see evidence of convergence to a limit.

We use the spectral data to study solutions of the heat equation and wave equation. Graphs of the heat Kernel in Figures 7.0-7.5 show the expected outward spreading, but it has to be remembered that the domain of the graphs is $OG$ as opposed to $POG$. Similarly we show the time evolution of a solution of the wave equation in Figures 8.0-8.17. "Motion pictures" of these solutions may be seen on the website [W].

In section 9 we examine geometric questions concerning the natural graph length metric of the cell graph approximations to $POG$. We address two basic questions: 1) what is the diameter $d_m$ of the level $m$ cell graph, and 2) what are the cardinalities $B_m(v, n)$ of balls in the level $m$ graph, where $v$ denotes the center point and $m$ the radius? For the first question we have a conjecture for the exact value $d_m = 2^{m-2} \cdot 7 + 1$ for $m \geq 2$. We are able to prove the estimates

$$2^{m-1} - 1 \leq d_m \leq 2^{m-1}(m+2) - 1$$

Note that the upper bound is $O(2^m m)$ rather than the conjectured $O(2^m)$. For the second question we are able to prove the upper bound $B_m(v, n) \leq 200n^3$ if $n < 2^{m-1}$. We have not been able to prove lower bound estimates. To obtain these estimates we developed intricate arguments concerning geodesics in the metric. Similar methods are used to study analogous questions concerning fractals obtained from the Sierpinski carpet by identifying boundaries using projective, Klein bottle and torus identifications in [GSS].

We all need to acknowledge two weaknesses in our work. The first is that, while it would be very valuable to extend our computations to higher levels, we do not see how to go beyond $m = 5$. It is the complexity of the identification algorithms described in section 4 that is the immediate barrier, even before we confront the memory space and computation time obstacles that exist at present, but can reasonably be expected to improve in the future. The other issue is the decision we made to use cell-graphs rather than vertex graphs in our approximations. This choice was motivated by the fear that the identification algorithms in section 4 would have to be even more intricate for the vertex graphs. Ideally it would have been better to implement both approaches to see if the computation outputs are similar. For now we must leave such improvements to the future.

**Acknowledgements:** Levente Szabo was supported by the Jane Oliver Swafford Scholarship through the University of North Carolina Asheville. Wing Hong Wong was supported by the 1978 Mathematics Alumnus Li Sze-lim Scholarship 2017/18 and the Yasumoto International Exchange Scholarship given through the Chinese University of Hong Kong. Finally we would like to thank the entire Cornell University staff and faculty for helping organize and facilitate the SPUR program and making this research possible.